\documentclass[12pt]{article}

\usepackage{arxiv}

\usepackage[utf8]{inputenc} 
\usepackage[T1]{fontenc}    
\usepackage{url}            
\usepackage{booktabs}       
\usepackage{amsfonts}       
\usepackage{nicefrac}       
\usepackage{microtype}      
\usepackage{lipsum}
\usepackage{graphicx}
\usepackage{amsmath}
\usepackage[style=apa,sortcites=true,sorting=nyt,backend=biber]{biblatex}
\graphicspath{ {./images/} }
\makeatletter
\g@addto@macro\bfseries{\boldmath}
\makeatother

\title{On the validity of bootstrap uncertainty estimates in the Mallows-Binomial model}

\author{
 Michael Pearce \\
  Department of Statistics\\
  University of Washington\\
  Seattle, WA \\
  \texttt{mpp790@uw.edu} \\
   \And
 Elena A. Erosheva \\
  Department of Statistics, School of Social Work, and the Center for \\Statistics and the Social Sciences\\
  University of Washington\\
  Seattle, WA \\
  \texttt{erosheva@uw.edu} \\
}

\begin{document}
\maketitle
\begin{abstract}
The Mallows-Binomial distribution is the first joint statistical model for rankings and ratings \parencite{pearce2022unified}. Because frequentist estimation of the model parameters and their uncertainty is challenging, it is natural to consider the nonparametric bootstrap. However, it is not clear that the nonparametric bootstrap is asymptotically valid in this setting. This is because the Mallows-Binomial model is parameterized by continuous quantities whose discrete order affects the likelihood. In this note, we demonstrate that bootstrap uncertainty of the maximum likelihood estimates in the Mallows-Binomial model are asymptotically valid.
\end{abstract}

\section{Introduction}

Rankings and ratings are two types of preference data, which are primarily studied separately in the literature. Rankings may be modeled via a variety of specialized distributions, such as the Mallows \parencite{Mallows1957}, Plackett-Luce \parencite{Plackett1975}, and Bradley-Terry \parencite{Bradley1952}. Ratings are not commonly modeled using statistical methods, but are instead studied using simple summary statistics such as the mean or median \parencite{lee2013bias,tay2020beyond,NIHPeerReview}. In recent years, a growing body of work suggests that modeling rankings and ratings jointly may improve the accuracy of preference modeling and preference aggregation \parencite{ovadia2004ratings,shah2018design,liu2022integrating}. The Mallows-Binomial is the first joint statistical model for rankings and ratings \parencite{pearce2022unified}.

\textcite{pearce2022unified} estimated the Mallows-Binomial model in the frequentist setting via the method of maximum likelihood. Similar to the standard Mallows model, estimation of the Mallows-Binomial is provably an NP-hard problem \parencite{Meila2012}. However, exact calculation of the MLE is possible when a relatively small number of objects are assessed via an exact tree-search algorithm based on A$^*$, and approximate algorithms can be used otherwise \parencite{pearce2022unified}. Still, analytic standard error results are challenging for the Mallows-Binomial and led the authors to propose using the nonparametric bootstrap for uncertainty quantification.

The bootstrap is a very general tool for uncertainty quantification that was first proposed in \textcite{efron1979}. The nonparametric bootstrap is used to estimate the inherent uncertainty of an estimator without making assumptions on its distributional form. Given $n$ i.i.d. observations, $X=(X_1,\dots,X_n)$, the nonparametric bootstrap is performed using the following steps:
\begin{enumerate}
    \item Re-sample $n$ observations with replacement from the original dataset
    $B$ times, for large $B$. Denote each bootstrap sample $X^b$, $b=1,\dots,B$.
    \item Estimate the unknown statistic(s) of interest, $\theta$, separately using each bootstrap sample. Denote the estimates from each bootstrap sample $\hat\theta^b$, $b=1,\dots,B$.
    \item Form an empirical distribution for $\hat\theta$ using the values $\hat\theta^b$, $b=1,\dots,B$.
\end{enumerate}
Quantiles from the empirical distribution of the unknown statistic(s) of interest may be used for the purpose of creating confidence regions.

Despite its wide applicability, the nonparametric bootstrap does not always yield asymptotically valid confidence intervals. A canonical example in which the bootstrap fails is the estimation of the unknown parameter $\theta$ given i.i.d. samples from a Uniform$(0,\theta)$ distribution. Here, the MLE of $\theta$ is the maximum order statistic, which has a limiting exponential distribution. The bootstrap empirical distribution, however, will not be able to replicate the asymptotic distribution given the fixed sample, in which the bootstrap estimates of $\theta$ will always be less than or equal to the full-sample maximum order statistic \parencite{bickel1981some}.  Another canonical example in which the bootstrap fails is in estimating the location parameter of a Cauchy distribution. In this setting, the MLE is the sample mean which is itself Cauchy distributed and therefore has infinite variance. As a result, the bootstrap estimator behaves poorly even given large samples \parencite{politis1998computer}.

The Mallows-Binomial likelihood has an unusual form that makes the asymptotic validity of bootstrap uncertainty for the MLE unclear. Specifically, the model is parameterized by continuous parameters whose discrete order impacts the likelihood. That is, the likelihood contains both continuous and discrete components; discontinuities may exist whenever the order of certain parameters change. Thus, frequentist estimation of the Mallows-Binomial model is both a continuous and discrete problem. In the absence of theoretical results regarding the asymptotic distribution of the maximum likelihood estimators, the validity of the nonparametric bootstrap is unclear. 

In this note, we demonstrate that the nonparametric bootstrap is an asymptotically valid method of uncertainty quantification for maximum likelihood estimates in the Mallows-Binomial model. The rest of the paper is organized as follows: In Section 2, we provide preliminaries regarding notation and a formal model statement. Main asymptotic results are provided in Section 3, followed by a conclusion in Section 4 that summarizes our work and suggests directions for future study.

\section{Preliminaries}

\subsection{Notation}

Suppose there exists a collection of $J$ objects which will be assessed by $I$ judges. Each judge rates each object using the integers $\{0,1,\dots,M\}$, where $M$ is a fixed and known maximum rating. Smaller ratings are better; higher ratings are worse. We let $X_{ij}$ represent the rating that judge $i$ assigns to object $j$. Additionally, each judge provides a ranking of their most-preferred objects. For simplicity, we assume each judge provides a ranking of all the objects. We let $\Pi_i$ denote judge $i$'s ranking.

Furthermore, suppose that each object $j\in\{1,\dots,J\}$ has a true underlying quality $p_{0j}\in(0,1)$. The vector of true underlying qualities is written as $p_0 = [p_{01} \ p_{02} \ \dots p_{0J}]^T$, such that $p_0\in(0,1)^J$. Given a vector $p\in(0,1)^J$, we let $\pi_p$ denote its order from least to greatest. For example, if $p = [0.3 \ 0.2 \ 0.1]^T$, then $\pi_p = 3\prec 2\prec 1$ (read as ``object 3 is preferred to object 2, which is preferred to object 1"). We call $\pi_{p_0}$ the true \textit{consensus ranking} of the objects. Additionally, we assume there is a true $\theta_0>0$, which represents the strength of population ranking consensus. $\theta_0$ is defined identically to that from a traditional Mallows $\phi$ model.

\subsection{Mallows-Binomial}

The Mallows-Binomial model is a joint distribution for rankings and ratings that is designed to capture the above situation. Let $(X,\Pi)_I$ denote an i.i.d. sample of size $I$ from a Mallows-Binomial($p_0,\theta_0$) distribution. Their joint likelihood can be written:
\begin{align}
    \mathcal{L}\Big((X,\Pi)_I|p_0,\theta_0\Big) &=\prod_{i=1}^I \Bigg(\frac{e^{-\theta_0  d(\pi_i,\pi_{p_0})}}{\psi(\theta_0)}\times\prod_{j=1}^J {M\choose x_{ij}}p_{0j}^{ x_{ij}}(1-p_{0j})^{M-x_{ij}}\Bigg)\label{eq:lik}\\
    &= \frac{e^{-\theta_0 \sum_{i=1}^I d(\pi_i,\pi_{p_0})}}{\psi(\theta_0)^I}\times\prod_{j=1}^J \Bigg(\prod_{i=1}^I {M\choose x_{ij}}\Bigg)p_{0j}^{\sum_{i=1}^I x_{ij}}(1-p_{0j})^{IM-\sum_{i=1}^I x_{ij}},
\end{align}
where $d(\cdot,\cdot)$ is the Kendall's $\tau$ distance between the two rankings and
\begin{equation}\label{eq:psi}
    \psi(\theta) = \prod_{j=1}^J \frac{1-e^{-j\theta}}{1-e^{-\theta}}
\end{equation}
is the normalizing function of a Mallows model. As can be seen from Equation \ref{eq:lik}, the likelihood of each observation corresponds to the product of a Mallows $\phi$ ranking distribution with $J$ Binomial rating distributions.

Next, we provide a preliminary expression for the MLE, $(\hat p,\hat\theta)$:

\begin{align}
    (\hat p,\hat \theta) &= \underset{p,\theta}{\arg\max} \Bigg[-\theta \sum_{i=1}^I d(\pi_i,\pi_p)-I\log\psi(\theta)\\
    & \hspace{65pt}+\sum_{j=1}^J (\sum_{i=1}^I x_{ij})\nonumber \log(p_j)+(IM-\sum_{i=1}^I x_{ij})\log(1-p_j)\Bigg]\\
    &= \underset{p,\theta}{\arg\max} \Bigg[-\theta \bar{D}(\pi,\pi_p)-\log\psi(\theta)+\sum_{j=1}^J\Big( \bar{x_{j}}\log(p_j)+(M-\bar{x_{j}})\log(1-p_j)\Big)\Bigg]\\
    &\equiv \underset{p,\theta}{\arg\max} \Bigg[f_{(X,\Pi)_I}(p,\theta)\Bigg],
\end{align}
where 
\begin{equation}
    \bar{D}(\pi,\pi_p)=I^{-1}\sum_{i} d(\pi_i,\pi_p),
\end{equation}
and
\begin{equation}
    \bar{x}_j = I^{-1}\sum_i x_{ij}.
\end{equation}

There is no closed-form solution for the MLE. However, computationally-efficient frequentist estimation is possible via a tree-search method based on the A$^*$ algorithm \parencite{pearce2022unified}.

In the remainder of this work, we will assume that in $p_0$, $p_{0j}\neq p_{0k}$ whenever $j\neq k$. Furthermore, we assume that each $p_{0j}\in(a,b)\subset [0,1]$, and that $\theta\in (c,d)\subset [0,\infty)$. Under these conditions, \textcite{pearce2022unified} proved that the maximum likelihood estimator $(\hat p,\hat\theta)$ is consistent for $(p_0,\theta_0)$ as the number of judges, $I$, grows to infinity. 

\section{Asymptotic Bootstrap Validity}

We would like to show that bootstrap uncertainty estimates are asymptotically valid for the MLE in the Mallows-Binomial model. A sufficient condition for bootstrap validity is local asymptotic normality of the MLE \parencite{hall2013bootstrap,bickel1981some}. As such, we show in the following subsections that the $(J+1)$-dimensional MLE $(\hat p,\hat\theta)$ is coordinate-wise, locally asymptotically normal.

\subsection{Local Asymptotic Normality of $\hat p_j$}\label{p}

We begin by considering each $\hat p_j$, $j=1,\dots,J$. Note that $\hat p_j$ is the solution to the following equation:
\begin{align}
    0 &= \frac{\partial}{\partial p_j} f_{(X,\Pi)_I}(p,\theta)\\
    &= -\theta \Big[\frac{\partial}{\partial p_j} \bar{D}(\pi,\pi_p)\Big] + \frac{\bar{x}_j}{p_j} -\frac{M-\bar{x}_j}{1-p_j}
\end{align}

A key challenge in calculating $\hat p_j$ is the derivative $\frac{\partial}{\partial p_j} \bar{D}(\pi,\pi_p)$, which is a function of the $J$-dimensional vector $p$. However, as long as each $p_k \neq p_j$ whenever $k\neq j$, then $\pi_p$ will remain constant in small perturbations around $p_j$ and the derivative in $p_j$ will thus be 0. 

Generalizing from $p_j$ to $p$, we require there to exist an $\epsilon$-ball around the J-dimensional vector $p$ such that $\text{Order}(p)=\pi_p$ remains constant for all $p'$ in the ball. In such cases, we have $\frac{\partial}{\partial p_j} \bar{D}(\pi,\pi_p)=0$ and thus $\hat p_j$ is defined by the standard Binomial MLE, $\bar{x}_j/M$, which is asymptotically normal as it is a function of the mean of i.i.d. random variables.

Specifically, this means that in a local region defined by the order of $\hat p$, we have the standard Binomial result,
\begin{align}
    \sqrt{I}\Big(\hat p_j - p_{0j}\Big) \overset{d}{\to} N\Big(0,\frac{p_{0j}(1-p_{0j})}{M}\Big),
\end{align}
which establishes the coordinate-wise local asymptotic normality of $\hat p_j$, $j=1,\dots,J$.

\subsection{Local Asymptotic Normality of $\hat\theta$}\label{theta}

We now show that $\hat\theta$ is coordinate-wise a locally asymptotically normal estimator of $\theta_0$. Note that $\hat\theta$ is the solution to the following equation:
\begin{align}
    0 &= \frac{\partial}{\partial\theta} f_{(X,\Pi)_I}(p,\theta)\\
    &= -\bar{D}(\pi,\pi_p) -\frac{\psi'(\theta)}{\psi(\theta)}\\
    \implies \bar{D}(\pi,\pi_p) &= -\frac{\psi'(\theta)}{\psi(\theta)}.
\end{align}

For simplicity, we define $\kappa(\theta) = -\psi'(\theta)/\psi(\theta)$. Thus, we have
\begin{align}
    \hat\theta &= \kappa^{-1}(\bar{D}(\pi,\pi_p)).
\end{align}
No simple expression for $\kappa^{-1}$ exists. However, it can be seen from Equation \ref{eq:psi} that when $J\geq2$ and for any $\theta>0$, $\psi$ is a continuous and positive function with a continuous and strictly negative first derivative and a continuous and strictly positive second derivative. Thus, $\kappa(\cdot)$ is a smooth and positive function. Furthermore, $\kappa$ is monotone decreasing \parencite{fligner1986distance}. As a result, its inverse $\kappa^{-1}(\cdot)$ is well-defined and so is $\hat\theta$ given $\bar{D}(\pi,\pi_p)$.

For reasons which will be made clear later, we also write out an expression for $\kappa(\theta)$:
\begin{align}
    \kappa(\theta) = -\frac{\psi'(\theta)}{\psi(\theta)}&= -\frac{\frac{(1-e^{-\theta})^J (\frac{\partial}{\partial\theta}\prod_{j=1}^J 1-e^{-j\theta})-(\prod_{j=1}^J 1-e^{-j\theta})Je^{-\theta}(1-e^{-\theta})^{J-1}}{(1-e^{-\theta})^{2J}}}{\frac{\prod_{j=1}^J 1-e^{-j\theta}}{(1-e^{-\theta})^J}}\\
    &= -\frac{(\frac{\partial}{\partial\theta}\prod_{j=1}^J 1-e^{-j\theta})-(\prod_{j=1}^J 1-e^{-j\theta})Je^{-\theta}/(1-e^{-\theta})}{\prod_{j=1}^J 1-e^{-j\theta}}\\
    &= \frac{Je^{-\theta}}{1-e^{-\theta}}-\frac{(\frac{\partial}{\partial\theta}\prod_{j=1}^J 1-e^{-j\theta})}{\prod_{j=1}^J 1-e^{-j\theta}}\\
    &= \frac{Je^{-\theta}}{1-e^{-\theta}}-\frac{\partial}{\partial\theta}\log(\prod_{j=1}^J 1-e^{-j\theta})\\
    &= \frac{Je^{-\theta}}{1-e^{-\theta}}-\sum_{j=1}^J \frac{je^{-j\theta}}{1-e^{-j\theta}}\label{eq:kappa}
\end{align}

Next, note that $\bar{D}(\pi,\pi)$ is a random variable given a fixed consensus ranking $\pi$, due to the randomness in the collection of rankings $\pi$. According to \textcite{fligner1986distance}, $\bar{D}(\pi,\pi_{p_0})$ is asymptotically normal with mean and variance depending on the true $\theta_0$. Specifically,
\begin{align}
    \mu_{\theta_0} &\equiv E_{\theta_0}[D(\pi_i,\pi_{p_0})] = \frac{Je^{-\theta_0}}{1-e^{-\theta_0}} - \sum_{j=1}^J \frac{je^{-j\theta_0}}{1-e^{-j\theta_0}}\label{eq:mu}\\
    \sigma^2_{\theta_0} &\equiv Var_{\theta_0}[D(\pi_i,\pi_{p_0})] = \frac{Je^{-\theta_0}}{(1-e^{-\theta_0})^2} - \sum_{j=1}^J \frac{j^2e^{-j\theta_0}}{(1-e^{-j\theta_0})^2}\\
    &\implies  \sqrt{I}(\bar{D}(\pi,\pi_{p_0}) -\mu_{\theta_0}) \overset{d}{\to} N(0,\sigma^2_{\theta_0})
\end{align}
Interestingly, we see from comparing Equations \ref{eq:kappa} and \ref{eq:mu} that $\mu_{\theta_{0}} = \kappa(\theta_0)$, which implies $\kappa^{-1}(\mu_{\theta_0})=\theta_0$. Since $\kappa^{-1\prime}$ is a real-valued function that does not equal 0, by the Delta method,
\begin{align}
    \sqrt{I}\Big(\kappa^{-1}(\bar{D}(\pi,\pi_{p_0})) -\kappa^{-1}(\mu_{\theta_0})\Big) &\overset{d}{\to} N\Big(0,\sigma^2_{\theta_0}[\kappa^{-1\prime}(\mu_{\theta_0})]^2\Big)\\
    \implies \sqrt{I}\Big(\hat\theta -\theta_0\Big) &\overset{d}{\to} N\Big(0,\sigma^2_{\theta_0}[\kappa^{-1\prime}(\mu_{\theta_0})]^2\Big)
\end{align}
Although the asymptotic variance cannot be written in closed-form expression, it is positive and finite. Therefore, in a local area of $\pi_{p_0}$, $\hat\theta$ is coordinate-wise an asymptotically normal estimator of $\theta$.

\subsection{Asymptotic Normality of $(\hat p,\hat\theta)$}

We showed in Sections \ref{p} and \ref{theta} that in a local neighborhood of the true consensus ranking $\pi_{p_0}$, the estimators $\hat p_j$, $j=1,\dots,J$ and $\hat\theta$ are coordinate-wise asymptotically normal. We add that the MLE $(\hat p,\hat\theta)$ is consistent as the number of observations, $I$, grows to infinity \parencite{pearce2022unified}. Thus, $\hat p \overset{p}{\to} p_0$ and $\hat \theta \overset{p}{\to} \theta_0$. As a result, as $I\to\infty$ the MLE $\pi_{\hat p}$ will be in a local neighborhood of the true $\pi_{p_0}$ with probability tending to 1, and the MLE $(\hat p,\hat\theta)$ will be coordinate-wise an asymptotically normal estimator of $(p_0,\theta_0)$ in that local neighborhood. This satisfies a sufficient condition for asymptotic bootstrap validity \parencite{hall2013bootstrap,bickel1981some}.

\section{Conclusion}

In this note, we demonstrate that bootstrap uncertainty estimates are asymptotically valid for the MLE in the Mallows-Binomial distribution. Our work only proves coordinate-wise, local asymptotic normality of the vector-valued MLE. This result technically guarantees only asymptotically accurate marginal coverage. As a result, the confidence intervals may suffer from either overcoverage or undercoverage when applied jointly. To understand why, we can think of the marginal confidence intervals jointly creating a confidence region that is a $(J+1)$-dimensional hypercube, as opposed to a $(J+1)$-dimensional confidence ellipse that could be created via a true joint analysis. Overcoverage may occur if the confidence hypercube contains as a subset the (theoretical) confidence ellipse. On the other hand, if each coordinate of the confidence hypercube is independent, joint coverage becomes a multiple testing problem and may result in undercoverage. That said, the present results do not preclude the possibility that bootstrap uncertainty estimates provide asymptotically correct coverage in the joint setting. We find no evidence to suggest asymptotically correct joint intervals are invalid. Further research may demonstrate proper coverage in the joint setting by establishing joint asymptotic normality of the $(J+1)$-dimensional MLE in a local neighborhood of $\pi_{p_0}$.

The present work does not address the computational challenges of frequentist estimation of the Mallows-Binomial model, which is an NP-hard problem \parencite{pearce2022unified,Meila2012}. As a result, forming an appropriate empirical distribution of the MLE $(\hat p,\hat\theta)$ via repeated estimation of bootstrap samples may be intractably slow. Bayesian estimation presents a natural alternative for estimating uncertainty in a unified framework, which may ultimately speed up the process of parameter estimation and inference and make moot the question of bootstrap validity.

\section*{Acknowledgements}

The authors would like to thank Yen-Chi Chen for his helpful feedback and advice while assembling this work. This work was funded by NSF Grant No. 2019901.

\printbibliography

\end{document}